\def \version {2019--1--24}

\documentclass[12pt]{article}

\usepackage{mathtools}
\usepackage{amsmath,amssymb,amsbsy}
\usepackage{graphicx}
\usepackage{multirow}

\usepackage{amsfonts}
\usepackage{amsthm}
\usepackage{amssymb}
\usepackage{graphicx, enumerate}
\usepackage{amsmath}
\usepackage{latexsym}
\usepackage{longtable}
\usepackage{tabularx}
\usepackage{amsmath}
\usepackage{amsfonts}
\usepackage{amsthm}
\usepackage{setspace}
\usepackage{graphicx}
\usepackage{float}
\usepackage{rotating}
\usepackage{tikz}
\usepackage{verbatim}

\newtheorem{thm}{Theorem}[section]
\newtheorem{lem}[thm]{Lemma}
\newtheorem{obs}[thm]{Observation}

\newtheorem{cor}[thm]{Corollary}

\newtheorem{conj}[thm]{Conjecture}
\newtheorem{rem}[thm]{Remark}

\makeatletter
\def\imod#1{\allowbreak\mkern10mu({\operator@font mod}\,\,#1)}
\makeatother
\def \vep {\varepsilon}
\def \cH {\mathcal{H}}

\def \bsk {\bigskip}

\newcommand{\zet}{\mathbb{Z}}

\begin{document}

\title{Realization of digraphs in Abelian groups and its consequences}

\author{Sylwia Cichacz$^{1,}$\footnote{This work was partially supported by the Faculty of Applied Mathematics AGH UST statutory tasks within subsidy of Ministry of Science and Higher Education.}~, ~Zsolt Tuza$^{2,3,}$\thanks{Research supported in part by
 the National Research, Development and Innovation Office -- NKFIH under the
  grant SNN 129364.}\\ ~~ \\
\normalsize $^1$AGH University of Science and Technology, \vspace{2mm} Poland\\
\normalsize $^2$Alfr\'{e}d R\'{e}nyi Institute of Mathematics, Hungarian Academy of Sciences\\
\normalsize Budapest, \vspace{2mm} Hungary\\
\normalsize $^3$Department of Computer Science and Systems Technology\\
\normalsize University of Pannonia, Veszpr\'em, Hungary}

%\author{Sylwia Cichacz\\
%\small Faculty of Applied Mathematics\\
%\small AGH University of Science and Technology\\
%\small Al. Mickiewicza 30, 30-059 Krak\'ow, Poland\\
%\small\tt cichacz@agh.edu.pl}
\date{\small Latest update on \version}
\maketitle
\begin{abstract}
Let  $\overrightarrow{G}$  be  a  directed  graph  with no component of order
 less than~$3$, and let $\Gamma$ be a finite Abelian group such that
  $|\Gamma|\geq 4|V(\overrightarrow{G})|$ or if $|V(\overrightarrow{G})|$ is large
 enough with respect to an arbitrarily fixed $\vep>0$ then
  $|\Gamma|\geq (1+\vep)|V(\overrightarrow{G})|$.
We  show that there  exists  an  injective  mapping  $\varphi$  from
$V(\overrightarrow{G})$  to  the  group $\Gamma$  such  that  $\sum_{x\in V(C)}\varphi(x) =  0$  for  every
connected  component  $C$  of  $\overrightarrow{G}$, where $0$ is the identity element of $\Gamma$.
Moreover we show some applications of this result to group distance magic labelings.
\end{abstract}

\vfil

 \section{Introduction}

Let  $\overrightarrow{G}=(V,A)$ be a directed graph. An arc $\overrightarrow{xy}$ is considered to be directed from $x$ to $y$, moreover $y$ is called the \emph{head} and $x$ is called the \emph{tail} of the arc. For a vertex $x$, the set of head endpoints adjacent to $x$ is denoted by $N^+(x)$, and the set of tail endpoints adjacent to $x$ is  denoted by $N^-(x)$.
 \bsk

Assume $\Gamma$ is an Abelian group of order $n$ with the operation denoted by~$+$.
  For convenience
we will write $ka$ to denote $a + a + \ldots + a$ where the element $a$ appears $k$ times, $-a$ to denote the inverse of $a$, and
we will use $a - b$ instead of $a+(-b)$.  Moreover, the notation $\sum_{a\in S}{a}$ will be used as a short form for $a_1+a_2+a_3+\dots$, where $a_1, a_2, a_3, \dots$ are all elements of the set $S$. The identity element of $\Gamma$ will be denoted by $0$. Recall that any group element $\iota\in\Gamma$ of order 2 (i.e., $\iota\neq 0$ and $2\iota=0$) is called an \emph{involution}.\bsk

Suppose that there exists a mapping $\psi$ from the arc set $E(\overrightarrow{G})$ of
$\overrightarrow{G}$ to an Abelian group $\Gamma$
such
that if we define a mapping $\varphi$ from the vertex set $V(\overrightarrow{G})$ of
$G$
to $\Gamma$
by
$$\varphi_{\psi}(x)=\sum_{y\in N^+(x)}\psi(yx)-\sum_{y\in N^-(x)}\psi(xy),\;\;\;(x\in V(G)),$$
then $\varphi_{\psi}$
is injective. In this situation, we say that $\overrightarrow{G}$ is \textit{realizable} in
$\Gamma$, and that the mapping $\psi$ is $\Gamma$-\textit{irregular}.\bsk% whereas $\varphi(x)$ a $\psi$-\textit{weighted degree of} $x$. \\

The corresponding problem in the case of simple graphs was considered
 in \cite{ref_AnhCic,ref_AnhCicMil,ref_AnhCicPrzy}.
  For  $\Gamma=(\zet_2)^m$ the problem was raised in \cite{Tuza}.
We easily see that if $\overrightarrow{G}$ is realizable in $(\zet_2)^m$,
 then every component of $\overrightarrow{G}$ has order at least 3
 (recall that we are assuming $\overrightarrow{G}$ has no isolated vertex).
The following results have been shown:

%Aigner and Triesch stated a conjecture that  $G$ is realizable in $\overrightarrow{G}$ if  $m = [\log_2 |V(G)|] + 1$. It was solved by Egawa
 
\begin{thm}[\cite{Egawa}] \label{Egawa} Let\/ $\overrightarrow{G}$ be a directed
 graph with no component of order less than\/ $3$.
Then\/ $\overrightarrow{G}$ is realizable in\/ $(\mathbb{Z}_2)^m$ if and only if\/
  $|V(\overrightarrow{G})|\leq 2^m$ and\/ $|V(\overrightarrow{G})|\neq 2^m-2$.
\end{thm}

\begin{thm}[\cite{Fukuchi}]  Let\/ $p$ be an odd prime and let\/
 $m\geq 1$ be an integer.
If\/~$\overrightarrow{G}$ is a directed graph without isolated vertices such that\/
 $|V(\overrightarrow{G})| \leq p^m$, then\/  $\overrightarrow{G}$ is realizable
  in\/ $(\zet_p)^m$.
\end{thm}

In this paper we will prove that a directed graph $\overrightarrow{G}$ with
 no component of order less than $3$ is realizable in any $\Gamma$ of order at
   least $|V(\overrightarrow{G})|$ such that either $\Gamma$ is of an odd order
   or $\Gamma$ contains exactly three involutions.
Moreover we will show that a directed  graph $\overrightarrow{G}$ with
 no component of order less than~$3$ is realizable in any $\Gamma$
  such that $|\Gamma|\geq 4|V(\overrightarrow{G})|$. Further, the coefficient 4 will be improved substantially for
 $|V(\overrightarrow{G})|$ large enough.
In the last section we will show some applications of this result.

\section{Characterizations and sufficient conditions}

A subset $S$ of $\Gamma$ is called a zero-sum subset if $\sum_{a\in S} a=0$.
It turns out that a realization of $\overrightarrow{G}$ in an Abelian group $\Gamma$
 is strongly connected with a zero-sum partition of $\Gamma$ \cite{ref_AigTri,Fukuchi}.
Using exactly the same arguments as in \cite{Fukuchi} for elementary Abelian groups
 we show the following for general Abelian groups.

\begin{thm}\label{zerosum} A directed graph\/ $\overrightarrow{G}$ with no isolated
 vertices  is realizable in\/ $\Gamma$ if and only if there exists an injective
mapping\/ $\varphi$ from $V(G)$ to $\Gamma$ such that\/
 $\sum_{x\in V\langle C)}\varphi(x)=0$ for every component\/ $C$ of\/ $G$.
\end{thm}

\textit{Proof.}
The  necessity  is  obvious.  To  prove  the  sufficiency,  let  $\varphi$  be  an  injective
mapping  from  $V(\overrightarrow{G})$  to  $\Gamma$  such  that   $\sum_{x\in V\langle C)}\varphi(x)=0$  for  every  connected  component  $C$  of  $\overrightarrow{G}$.

Let  $C$ be  a  connected  component  of  $\overrightarrow{G}$.  It  suffices  to  show  that  there
exists  a  mapping  $\psi$ from  $E(C)$  to  $\Gamma$  satisfying:
$$\varphi(x)=\sum_{y\in N^+(x)}\psi(yx)-\sum_{y\in N^-(x)}\psi(xy),\;\;\;(x\in V(G)).$$
Now we will construct a spanning tree of $C$. Let $V(C)  =\{  x_1  ,\ldots,x_k  \}$  ($k  =  |V(C)|$)
so  that  for  each  $2 \le i  <  k$,  there  exists  exactly one  arc $e_i$  between  $\{x_l,\ldots ,  x_{i-1}\}$  and  $x_i$.  For  each  $2 \le i  <  k$,  define  a  subdigraph  $C_i$  of  $C$ by  setting  $V(C_i)  =  V(C)$  and  $E(C_i)  =\{  e_{i+l},\ldots,e_k \}$. Let $\psi(e_k)=\varphi(x_k)$.  We  define  $\psi$  backward  inductively  by
$$\psi(e_i)=\left\{\begin{array}{l}
\varphi(x_i)  -\sum_{y\in N^+_{C_i}(x_i)}\psi(x_iy)  +\sum_{y\in N^-_{C_i}(x_i)} \psi(yx_i),\; \mathrm{if}\; x_i\; \mathrm{is\;the\;tail\;of\;} e_i,\\
-\varphi(x_i)  +\sum_{y\in N^+_{C_i}(x_i)}\psi(x_iy)  -\sum_{y\in N^-_{C_i}(x_i)} \psi(yx_i),\; \mathrm{if}\; x_i\; \mathrm{is\;the\;head\; of}\; e_i.
 \end{array}\right.$$
Finally,  let  $\psi(e)  =  0$  for  all  $e  \in  E(C)  \setminus \{e_2 , \ldots  ,  e_m  \}$.  Then  the  resulting  mapping  $\psi$
has  the  desired  property.~\qed\bigskip

The following result is known.

\begin{thm}[\cite{KLR, Zeng}]\label{Zeng} Let\/ $\Gamma$ have order\/ $n$.
For every partition\/ $n-1 = r_1 + r_2 + \ldots + r_t$ of\/ $n-1$ with\/
 {$r_i \geq 2$} for\/ $1 \leq i \leq t$ and for any possible positive integer\/ $t$,
  there is a partition of\/ $\Gamma-\{0\}$ into pairwise disjoint subsets\/
   $A_1, A_2,\ldots , A_t$ such that\/ $|A_i| = r_i$ and\/ $\sum_{a\in A_i}a = 0$
    for all\/ $1 \leq i \leq t$ if and only if either\/ $\Gamma$ is of an odd order
     or\/ $\Gamma$ contains exactly three involutions.
\end{thm}

By Theorems \ref{zerosum} and \ref{Zeng} we obtain the following immediately.

\begin{thm}A directed graph\/ $\overrightarrow{G}$ with no component of order
 less than\/ $3$ is realizable in any\/ $\Gamma$ of order at
 least\/ $|V(\overrightarrow{G})|$ such that either\/ $\Gamma$ is of an
  odd order or\/ $\Gamma$ contains exactly three involutions.
\end{thm}

Before we proceed to groups having more than three involutions,
 we need some lemmas.
For the sake of simplicity, for any element $a\in\Gamma$, we are going to use the notation $a/2$ for an arbitrarily chosen element $b\in\Gamma$ satisfying $2b=a$.  Let $S_{a/2}=\{b \in \Gamma: 2b =a \}$. %It was proved:
%\begin{myobservation}[\cite{ref_AnhCic}]\label{pierwiastek}
%If $\Gamma$ is an Abelian group of odd order, then $|S_{a/2}|=1$ for any element $a\in\Gamma$.
%\end{myobservation}

\begin{obs}\label{pierwiastek2}
If\/ $\Gamma$ is an Abelian group of even order\/ $n$, then\/ $|S_{a/2}|\leq n/2$
 for any element\/ $a\in\Gamma$, $a\neq 0$.
\end{obs}
\textit{Proof.}
If for some $a\neq 0$ there exist $g_1,g_2\in\Gamma$, $g_1\neq g_2$ such that $2g_1=a$ and $2g_2=a$ then it follows that $2(g_1-g_2)=0$ and consequently $g_1-g_2$ is an involution. Since $2g_1=a\neq 0$,
   the number of involutions in $\Gamma$ is
    less than $|\Gamma|/2$.\qed

\begin{lem}Let\/ $\overrightarrow{G}$ be a directed  graph with no component of
 order less than\/~$3$, and let\/ $\Gamma$ be a finite Abelian group such that\/
  $|\Gamma|\geq 4|V(\overrightarrow{G})|$.
There exists a\/ $\Gamma$-irregular labeling\/ $\psi$ of\/ $\overrightarrow{G}$
 such that\/ ${\psi}(e)\neq 0$ for every\/ $e\in E(\overrightarrow{G})$, and\/
  $\varphi_{\psi}(x)\neq 0$ for every\/ $x\in V(\overrightarrow{G})$.
\end{lem}
\textit{Proof.} The proof follows by induction on the number of arcs.

Suppose first that $\overrightarrow{G}$ is a path %Let us start with smallest possible graph $G$, i.e.
$\overrightarrow{P}_{3}$ with vertices, say, $u$, $v$ and $w$ and arcs $e_1$ and $e_2$. With no loss of generality we can assume that $e_1\cap e_2=v$. Let $\Gamma$ be an arbitrary Abelian group of order at least $12$.  Set an element $a\neq 0$ in such a way that $\varphi_{\psi}(u)=a$ (namely, $\psi(vu)=a$ and $\psi(uv)=-a$). Now, choose any $b\not\in\{0,a,-a,-2a\}$ and $b\not\in S_{-a/2}$. The number of forbidden values is at most $4+|\Gamma|/2<|\Gamma|$. Set now the element $b$ in such a way that $\varphi_{\psi}(v)=-a-b$. Both arc labels are different from $0$, and so are the vertex weighted degrees, since $\varphi_{\psi}(u)=a$, $\varphi_{\psi}(v)=-a-b$ and $\varphi_{\psi}(w)=b$. It is also obvious that the weighted degrees are three distinct elements of $\Gamma$.

Now let $\overrightarrow{G}$ be arbitrary directed graph of order $n$ with at least $3$ edges, having no component of order less than $3$, and let $\Gamma$ be any Abelian group of order at least $4n$. In the induction step we can assume that for every proper subgraph $\overrightarrow{H}$ of $\overrightarrow{G}$ having no component of order less than $3$ and for every Abelian group $\Gamma'$ of order at least $4|\overrightarrow{H}|$, there is a $\Gamma'$-irregular labeling ${\psi}_H$ of $\overrightarrow{H}$ in which no edge has label $0$ and $\varphi_{\psi_H}(x)\neq 0$ for every $x\in V(\overrightarrow{H})$. In particular, there is such labeling of $\overrightarrow{H}$ with $\Gamma^\prime=\Gamma$, since $|\Gamma|\geq 4n\geq 4|V(\overrightarrow{H})|$. We will extend ${\psi}_H$ to the labeling ${\psi}$ of $\overrightarrow{G}$, having the same properties.

We choose $\overrightarrow{H}$ in one of the following ways. If there is a component $C\cong \overrightarrow{P}_{3}$ of $\overrightarrow{G}$, then $\overrightarrow{H}=\overrightarrow{G}-C$. Otherwise, if there is a component $C$ and an edge $e\in E(C)$ not being a bridge in $C$, then $\overrightarrow{H}=\overrightarrow{G}-e$. Finally, if $\overrightarrow{G}$ is a forest with each component of order at least $4$, then choose any leaf edge $e$ of any component and let $\overrightarrow{H}=\overrightarrow{G}-e$.

Let us consider the first case. Assume that $\overrightarrow{G}=\overrightarrow{H}\cup \overrightarrow{H}^\prime$, where $V(H^\prime)=\{u,v,w\}$ and $E(H^\prime)=\{e_1,e_2\}$ such that $e_1\cap e_2=u$. Let ${\psi}_H$ be a  $\Gamma$-irregular labeling of $\overrightarrow{H}$ fulfilling the desired non-zero properties, existing by the induction hypothesis. Let ${\psi}(e)={\psi}_H(e)$ for $e\in E(\overrightarrow{H})$. Now choose any element of $a\in \Gamma$ such that $a\neq 0$ and $a\neq \varphi_{\psi_H}(x)$ for $x\in V(\overrightarrow{H})$ and set the label $a$ on the edge $e_1$ such that $\varphi_{\psi}(v)=a$. Such $a$ can be chosen, as only $n-3$ vertex weighted degrees have been assigned so far and $|\Gamma|>n-2$. Now choose $b\in\Gamma$ such that $b\not\in \{0,a,-a,-2a\}$, $b\not \in S_{-a/2}$, $ b \not \in \{\varphi_{\psi_H}(x), -w(x)+a\}$ for $x\in V(H)$ and set $b$ on the arc $e_2$ such that $\varphi_{\psi}(w)=b$. The number of forbidden elements is at most $4+|\Gamma|/2+2(n-3)=2n-2+|\Gamma|/2<|\Gamma|$, so we can choose such $b$. Obviously, the two new edge labels are not $0$ and neither are the three new weighted degrees $\varphi_{\psi}(v)=a$, $\varphi_{\psi}(w)=b$ and $\varphi_{\psi}(u)=-a-b$. Also, the three new weighted degrees are pairwise distinct and not equal to any $\varphi_{\psi}(x)$, where $x\in V(\overrightarrow{H})$.

In the second case, let $\overrightarrow{H}=\overrightarrow{G}-e$ and let ${\psi}_H$ be a  $\Gamma$-irregular labeling of $\overrightarrow{H}$ fulfilling the desired non-zero properties, existing by the induction hypothesis. Now let ${\psi}(y)={\psi}_H(y)$ for $y\in E(\overrightarrow{H})$. Let us denote the tail of $e$ by $u$ and the head by $v$. Choose an element $a\in\Gamma$ such that $a\not\in \{0,\varphi_{\psi_H}(u),-\varphi_{\psi_H}(v)\}$, $a\neq \varphi_{\psi_H}(u)-\varphi_{\psi_H}(x)$ for $x\in V(\overrightarrow{G})\setminus \{u,v\}$ and $a\neq \varphi_{\psi_H}(x)-\varphi_{\psi_H}(v)$ for $x\in V(\overrightarrow{G})\setminus \{u,v\}$, and $a\not\in S_{(\varphi_{\psi_H}(u)-\varphi_{\psi_H}(v))/2}$. Set ${\psi}(uv)=a$. The number of forbidden values is at most $3+2(n-2)+|\Gamma|/2<|\Gamma|$, so we can always choose such $a$. Note that two adjusted weighted degrees remain distinct and because of the way that $a$ was chosen, they are different from any weighted degree $\varphi_{\psi}(x)$ for $x\in V(\overrightarrow{G})\setminus \{u,v\}$. This means that $\psi$ has the desired property.

Finally, consider the third case. Assume that the ends of $e$ are $u$ and $v$, where $u$ is the pendant vertex. Having a  $\Gamma$-irregular labeling $\psi_H$ of $\overrightarrow{H}$ fulfilling the desired non-zero properties, we set $\psi(y)=\psi_H(y)$ for $y\in E(\overrightarrow{H})$. Then we choose $a\in \Gamma$ such that $a\not\in \{0,-\varphi_{\psi_H}(v)\}$,  $a\neq \varphi_{\psi_H}(x)$ for $x\in V(\overrightarrow{G})\setminus \{u,v\}$ and $a\neq \varphi_{\psi_H}(v)-\varphi_{\psi_H}(x)$ for $x\in V(\overrightarrow{G})\setminus \{u,v\}$, and $a\not\in S_{\varphi_{\psi_H}(u)/2}$. There are at most $2+2(n-2)+|\Gamma|/2<|\Gamma|$ forbidden values, so we can choose such $a$. Put label $a$ on the edge $e$ such that $\varphi_{\psi}(u)=a$. The adjusted weighted degree $\varphi_{\psi}(v)$ and the new weighted degree $\varphi_{\psi}(u)$ are distinct and different from any of the weighted degrees $\varphi_{\psi}(x)$ for $x\in V(G)\setminus \{u,v\}$, so also in this case $\overrightarrow{G}$ has the labeling $\psi$ with the desired property. This completes the proof. \qed

\bigskip

The above lemma implies the following.
\begin{thm}\label{glowne}A directed  graph\/ $\overrightarrow{G}$ with no component
 of order less than\/ $3$ is realizable in any\/ $\Gamma$ such that\/
  $|\Gamma|\geq 4|V(\overrightarrow{G})|$.
\end{thm}

\section{Asymptotic result}

Our goal here is to prove that if $|\Gamma|$ gets large, then it is possible to
 strengthen Theorem \ref{glowne} considerably, by replacing the
 multiplicative constant 4 in the condition $|\Gamma|\geq 4|V(\overrightarrow{G})|$
 with $(1+o(1))$, and also omitting the assumption that $\Gamma$ has more than one involution.
In the proof we shall apply the following corollary of Theorem~1.1
 from \cite{FraRod}.

\begin{lem}   \label{l:FR}
For every fixed\/ $\vep>0$ there exists an\/ $n_0=n_0(\vep)$
 with the following properties.
If\/ $\cH$ is a 3-uniform regular hypergraph with\/ $n>n_0$ vertices such that
 the degree of regularity is at least\/ $n/3$,
%  {{\komm (!!!) \footnote{SC - $n/3$? ---
%  ZS: I think you are right, maybe I was too much on the safe side.
%  The worst case seems to be $|R|=2|I|+2$, and then the degree, for which we apply
%  the lemma, is $(|R|-|I|-2)/2=|I|$. If $|I|\geq 3$ then indeed $|I|>(2|I|+2)/3$
%  and we may write n/3 instead of n/4 (if there is no error in this computation).
%  If R is larger and I is smaller, then (R-I-2)/2R = 1/2 - (I+2)/2R becomes larger,
%  hence remains above 1/3. Actually 3/8 seems to be good enough, but 1/3 is
%  simpler to state. And if we change n/4 to n/3, then also the same has to be
%  modified in the paragraph after the lemma.--- SC - I checked the inequalities
% and they seems to be OK} }},
  and each vertex pair is contained in at most two hyperedges,
 then\/ $\cH$ contains at least\/ $(1/3-\vep/3)n$ pairwise disjoint hyperedges.
Moreover if\/ $\cH$ is a 4-uniform hypergraph with\/ $n>n_0$ vertices, such that
 the vertex degrees are nearly equal and at least\/
 $n/4$, and each vertex pair is contained in at most three hyperedges,
 then\/ $\cH$ contains at least\/ $(1/4-\vep/4)n$ pairwise disjoint hyperedges.
\end{lem}

In fact the degree condition $n/3$ (and also $n/4$ in the 4-uniform case)
 can be replaced with $cn$ with any constant
 $c>0$, but we shall not need this stronger version of the lemma
  to allow very small degrees.

As another tool, we will use the following corollary of Theorem~\ref{Egawa}.

\begin{cor}[\cite{Egawa}]   \label{l:Egava2}
Let\/ $p\geq 2$ be an integer, and let\/ $q_3,$ $q_4,$ $q_5$ be nonnegative integers
 such that\/ $3q_3+4q_4+5q_5\leq 2^{p}$ and\/ $3q_3+4q_4+5q_5\neq 2^{p}-2$.
Then there exists a family\/ $Z=\{S_1,\dots,S_{q_3+q_4+q_5}\}$ of\/ $q_3+q_4+q_5$
 mutually disjoint zero-sum subsets of\/ $(\zet_2)^p$ such that\/
 $|S_i|=3$ for all\/ $1\le i\le q_3$,
 $|S_i|=4$ for all\/ $q_3+1\le i\le q_3+q_4$, and\/
 $|S_i|=5$ for all\/ $q_3+q_4+1\le i\le q_3+q_4+q_5$.
\end{cor}

The main result of this section is the following.

\begin{thm}   \label{t:pairs2}
Let\/ $\vep>0$ be fixed and assume that\/ $n$ is sufficiently large
 with respect to\/ $\vep$.
Let\/ $\Gamma\ne (\zet_2)^m$ be of order\/ $n$, and consider any integers\/
 $r_1,r_2,\dots,r_t$ with\/ $n > (1+\vep) (r_1 + r_2 + \ldots + r_t)$
  and\/ {$r_i \geq 2$} for all\/ $1 \leq i \leq t$.
If\/ $r_i = 2$ holds for at most\/ $n/4$ terms\/ $r_i$, then
 there exist pairwise disjoint subsets\/ $A_1, A_2,\ldots , A_t$ in\/ $\Gamma-\{0\}$
 such that\/ $|A_i| = r_i$ and\/ $\sum_{a\in A_i}a = 0$ for\/ $1 \leq i \leq t$.
\end{thm}

\newcounter{pnt}
\newcommand{\pont}%[1]
{\medskip \noindent ${%#1
\stepcounter{pnt} \arabic{pnt}}^{\circ}$\quad }
\textit{Proof.}
In order to make the structure of the argument more transparent,
 we split it into several parts.

\pont Assume $n>(2/\vep)\cdot n_0(\vep/2)$, where the
 function $n_0$ is from Lemma \ref{l:FR}.
We denote by $I$ the set of involutions in $\Gamma$,
 and write $R$ for the set of the other nonzero elements, i.e.\
  $R = \Gamma \setminus (I\cup\{0\})$.
We shall distinguish between the elements $\iota_1,\iota_2,\dots$ of $I$
 with subscripts, and wrtite $a,b,c,\dots$ for the elements of $R$.
A generic element may simply be denoted by $\iota\in I$ or $a\in R$.
 
 Since $\Gamma\ne (\zet_2)^m$, we have $|R|\geq n /2$.
If $|I|\leq \vep n /2$, then we omit $I$,
 and continue work in $R$ alone.
(This simplification also involves the elimination of the involution
 in case if $\Gamma$ has just one; then of course the elements of $R$
 sum up to 0.)
Otherwise we have both $|I|$ and $|R|$ larger than $n_0(\vep/2)$.
Below we describe the procedure for this more general case.

\pont If there are terms $r_i$ larger than 6, we modify the sequence by
 splitting each large term into a combination of terms 3 and 4.
Once the new sequence admits suitable zero-sum subsets, a solution for
 the original sequence follows immediately.
Note that this step does not create any new $r_i=2$, i.e.\ the condition
 on the number of terms 2 does not get violated.
Consequently we may assume $r_i\in \{2,3,4,5\}$ for all $1\le i\le t$.

In order to make further simplification, we state and prove the theorem
 in the following stronger form:
  \begin{itemize}
    \item[$(\star)$] The required disjoint zero-sum subsets $A_i$ exist also
     under the weaker assumption that the number of $r_i=2$ terms
     is at most $|R|/2$.
  \end{itemize}
Note that the bound $|R|/2$ is absolutely tight for every $\Gamma$
 because a zero-sum pair necessarily is of the type $(a,-a)$.

Now, if there is an $r_i=5$, we may split it into $2+3$, unless there are
 exactly $|R|/2$ terms $r_i=2$.
Similarly, if there is an $r_i=4$, we may split it into $2+2$, unless the
 number of $r_i=2$ terms is $|R|/2$ or $|R|/2-1$.
In this way the family of sequences $r_1,\dots,r_t$ to be studied is
 reduced to the following two cases:
  \begin{enumerate}
    \item there are exactly $|R|/2$ or $|R|/2-1$ terms $r_i=2$, and
     all the other terms are 3, 4, or 5; or
    \item we have $r_i\in \{2,3\}$ for all $1\le i\le t$.
  \end{enumerate}

\pont In case of 1., we can obtain the following much stronger result:
  \begin{itemize}
    \item[$(\star\star)$] If the number of $r_i=2$ terms is $|R|/2$ or
     $|R|/2-1$, then the required disjoint zero-sum subsets $A_i$ exist
     whenever $|\Gamma|\ge r_1+\dots+r_t+5$.
  \end{itemize}
Indeed, if $|I|=1$ then we lose at most two elements
 from $R$ and the only one element of $I$.
Otherwise $\Gamma$ has at least three involutions and we may apply
 Corollary \ref{l:Egava2}.
Namely, the 3-, 4-, and 5-terms can surely be assigned to suitable subsets
 $A_i$ of $I$ if their sum is at most $|I|-2$; and
 creating the $(a,-a)$ pairs for $a\in R$ we lose at most two elements
 from $R$ (and the 0-element of~$\Gamma$).

\pont In order to handle the case 2., we create an auxiliary set $R^*$ whose
 elements represent the inverse pairs of $R$, i.e.\ each $a^*\in R^*$
  stands for $(a,-a)$; hence $|R^*|=|R|/2$.
Note that each inverse-free subset of $R$ defines a unique subset of $R^*$
 with the same cardinality, in the natural way, while a $k$-element subset of
 $R^*$ may arise from $2^k$ distinct subsets of $R$.
One should be warned, however, that $R^*$ does not inherit the group
 structure of $R$.
Indeed, if $a+b=c$ holds inside $R$, then $a+(-b)\neq(-c)$; that is,
 $b^*=(-b)^*$, but $(a+b)^*\neq (a-b)^*$.

\pont Using Corollary \ref{l:Egava2} again, inside $I$ we define a large
 family $T_I$ of pairwise disjoint triples which together nearly cover $I$,
  such that the sum of the three elements in each triple equals~0.
If $|I|$ is a multiple of 3, we can partition $I$ into 3-element zero-sum
 subsets.
Otherwise, if $|I|=2^p-1=3q+1$, we find $q-1$ triples and one quadruple,
 each of whose elements sum up to zero.
Thus $T_I$ covers all but at most four elements of $I$.

\pont An analogous set $T_R$ of pairwise disjoint triples which nearly
 cover $R$ is more complicated to construct, because we put \emph{two}
 requirements instead of just one: if a triple $\{a,b,c\}$ is in $T_R$, then
  \begin{itemize}
   \item $a+b+c=0$; and
   \item the inverse triple $\{-a,-b,-c\}$ also belongs to $T_R$.
  \end{itemize}

We first construct an edge-labeled complete graph whose vertex set is $R$,
 and each edge $ab\in \binom{R}{2}$ gets the label
  $\lambda(a,b) := (-a)+(-b)$.
Hence $\lambda(-a,a)=0$ for all $a$ by definition,
 we shall disregard these edges.
At each $a\in R$ precisely $|I|$ edges are labeled from $I$, and
 consequently $|R|-|I|-2$ edges are labeled from $R$.
The strategy depends on whether the former or the latter is larger.

If $|I|<|R|/2$, or equivalently $|R|\geq 2|I|+2$,
 we keep the edges labeled from $R$.
It means that for all such edges we have
 $T_{a,b} := (a,b,\lambda(a,b)) \in \binom{R}{3}$, moreover every
 $T_{a,b}$ is a zero-sum triple.
This gives rise to the 3-uniform hypergraph, say $H_R$, whose hyperedges are
 the triples $T_{a,b}$, each pair of vertices belonging to
 either 0 or 1 hyperedge, and therefore each vertex being incident
  with exactly $(|R|-|I|-2)/2$
   hyperedges, hence the vertex degrees satisfy the inequality
   $\frac{(|R|-|I|-2)/2}{|R|} \geq
    1/2 - \frac{|I|+2}{4|I|+4} \geq 1/2 - 5/16 = 3/16$.

Note that if $T_{a,b}\in H_R$ then also $T_{-a,-b}\in H_R$
 (and of course vice versa), and they yield the same triple
 $T_{(a,b)^*}=T_{(-a,-b)^*}$ in the corresponding system $H_{R^*}$ over $R^*$.
We observe that inside the 6-tuple $T_{a,b} \cup T_{-a,-b}$
 there do not exist any further zero-sum triples.
Indeed, a third such triple should contain at least one element from
 each of $T_{a,b}$ and $T_{-a,-b}$, hence it would
 be of the form $T_{a,-b}$ or alike.
But then we would have $b-a\in \{a+b, -a-b\}$, from where $a+a=0$
 or $b+b=0$ would follow, contrary to the assumption $a,b\notin I$.
(This is in agreement with the comment above that $R^*$ does not inherit
 the group structure of $R$.)
It follows that the number of triples incident with a vertex in $H_{R^*}$
 is exactly the same as that in $H_R$.
In particular, $H_{R^*}$ is regular of a degree at least $\frac{3}{8}\,|R^*|$.
Moreover the maximum number of triples containing a pair $a^*,b^*$
 increases from 1 to 2, but not more.
Therefore Lemma~\ref{l:FR} can be applied and we obtain
 $(1/3-\vep/6)\cdot|R^*|$ pairwise disjoint triples in $R^*$.
Each of those triples $(a^*,b^*,c^*)$ originates from a triple $(a,b,c)$
 with $a+b+c=0$, hence it generates $(a,b,c)$ and $(-a,-b,-c)$
 inside $R$.
We denote this collection of disjoint zero-sum triples by $T_R$.
They together cover $(1-\vep/2)\cdot|R|$ elements of $R$
  in the case of $|I|<|R|/2$.

Otherwise, if $|I|\geq |R|/2$, note first that $|R|=|I|+1=|\Gamma|/2$ must hold,
 because $|R|$ is divisible by $|I|+1$ in every $\Gamma$.
Indeed, $a+\iota\in R$ holds for all
 $a\in R$ and $\iota\in I$.
Moreover, if $b-a=\iota_1$ and $c-b=\iota_2$ then $c-a=\iota_1+\iota_2\in I$
 hence the reflexive closure of the relation `\,$b-a\in I$\,' is an equivalence
 relation over $R$ and each of its equivalence classes contains exactly
 $|I|+1$ elements.

Consequently, the set $\{a+\iota\mid \iota\in I\}$ is the same as $R\setminus\{a\}$,
 therefore we have $-a=a+\iota_0$ for some $\iota_0\in I$
  (and of course $a=(-a)+\iota_0$ also).
 We observe that $-b=b+\iota_0$ holds for all $b\in R$.
Indeed, if $b=a+\iota_1$, then
 $-b=\iota_1-a=\iota_1+(-a)=\iota_1+a+\iota_0=b+\iota_0$.
(In case of larger $R$, this property would be guaranteed inside each equivalence
 class only.)

Let us put $x:=\lceil |R|/6 \rceil$, and take $x$ triples from $T_I$
 constructed above, such that none of them contains $\iota_0\in I$.
This selection can be done because $|T_I|\geq (|I|-4)/3 = (|R|-5)/3\ge |R|/6+1$
 whenever $\Gamma$ is not too small.
Recall that each member of $T_I$ is of the form $(\iota_1,\iota_2,\iota_3)$
 with $\iota_1+\iota_2+\iota_3=0$.
We represent the $x$ selected triples with new \emph{vertices} $u_1,u_2,\dots,u_x$,
 and define a nearly regular 4-uniform hypergraph on the vertex set
 $Q^*:=R^* \cup \{u_1,u_2,\dots,u_x\}$.
The hyperedges are of the form $(a^*,b^*,c^*,u_j)$, where $a$
 is any element, $b=a+\iota_1$, $c=b+\iota_2$ (hence $a=c+\iota_3$).
We do this for all $a\in R$ and all permutations of $\iota_1,\iota_2,\iota_3$
 in all of the first $x$ triples.

It is important to note that $a^*,b^*,c^*$ are three distinct elements
 because the triple of $T_I$ containing $\iota_0$ has not been selected.
For the same reason, for any fixed permutation $(\iota_1,\iota_2,\iota_3)$
 of any selected triple, the elements $a$ and $-a$ yield the same quadruple;
 that is, disregarding the representing new vertex $u_j$, we have
 $\{a^*,(a+\iota_1)^*,(a+\iota_1+\iota_2)^*\} =
   \{(-a)^*,(-a+\iota_1)^*,(-a+\iota_1+\iota_2)^*\}$.
Moreover, since any two of $\iota_1,\iota_2,\iota_3$ sum to the third,
 the quadruples of the form $\{a^*,(a+\iota_1)^*,(a+\iota_2)^*,(a+\iota_3)^*\}$
 partition $R^*$; this holds for each $u_j$.
Inside each such quadruple, three of the four 3-element subsets contain $a^*$.
  It follows that there are exactly $|R^*|=|R|/2$ hyperedges incident
 with any $u_j$, and the degree of an $a^*$ equals $3x\approx |R|/2$.
Thus, the conditions of Lemma \ref{l:FR} are satisfied, and
 there is a large packing of 4-element hyperedges $(a^*,b^*,c^*,u_j)$
 covering all but at most $(\vep/2)\cdot|Q^*|$ vertices of $Q^*$.

For every $(a^*,b^*,c^*,u_j)$ and its corresponding $(\iota_1,\iota_2,\iota_3)$
 we create the triples\footnote{Their inverse triples
  $(-a,b,\iota_1)$, $(-b,c,\iota_2)$, $(-c,a,\iota_3)$ would be equally fine.}
  $(a,-b,\iota_1)$, $(b,-c,\iota_2)$, $(c,-a,\iota_3)$.
All these three are zero-sum triples, and they partition the 9-tuple
 $\{a,b,c,-a,-b,-c,\iota_1,\iota_2,\iota_3\}$.
Observe further that no triangle from $T_I$ is used more than once in the
 construction; this is ensured by the presence of vertices $u_j$.
In this way we obtain $T_R$ if $|I|\geq|R|/2$.

\pont In case 2., we have $r_i\in\{2,3\}$ for all $1\le i\le t$.
 Let $m_k$ denote the number of terms $r_i=k$ for $k=2,3$.

If $m_3\le |T_I|$, we simply take any $m_3$ triples from $T_I$, and
 choose $m_2\le|R|/2$ pairs $(a,-a)$ inside $R$.
Otherwise two different situations may occur,
 depending on whether $|I|<|R|/2$ or not.
In both cases we assume that $m_3>|T_I|$ holds.

If $|I|<|R|/2$, we take the triples of $T_I$, moreover
 $2\cdot\lceil (m_3-|T_I|)/2 \rceil$
 triples from $T_R$ in such a way that if a triple $(a,b,c)$ is selected,
 then we also select $(-a,-b,-c)$.
This may yield one more triple than what we need, which we shall forget
 at the very end; but currently it is kept, in order to ensure that
 the rest of $R$ consist of inverse pairs.

If $|I|\geq|R|/2$, we again start with the triples of $T_I$, but then replace
 $\lceil (m_3-|T_I|)/2 \rceil$ of them with three triples from $T_R$ each.
This can be done by choosing $\lceil (m_3-|T_I|)/2 \rceil$ from the first
 $x$ members of $T_I$, and replacing them with the triples covering
 $\{a,b,c,-a,-b,-c,\iota_1,\iota_2,\iota_3\}$ as constructed above.

In either case, since $I$ is covered with the exception of at most four
 elements, there remains enough room for selecting the $m_2$ pairs
 $(a,-a)$ in the part of $R$ which is not covered by the selected triples.
% 
% \pont So far we have left at most $\vep n/2$ elements of $I$ uncovered.
% Thus, for $n > (1+\vep) (r_1 + r_2 + \ldots + r_t)$ we can select
%  the required $(1/4-\vep)n$ or fewer zero-sum pairs $(a,-a)$
%  in the uncovered part of $R$.
% 
% Summarizing, the number of uncovered or omitted elements (omission
%  occurs if $I$ or $R$ is small) is at most
%  $\vep n/2$ in each of $I$ and $R$, therefore we have enough room for
%  the zero-sum subsets $A_1,A_2,\dots,A_t$.
This completes the proof of the theorem.~\qed

\begin{cor}
There exists a function\/ $h_0:\mathbb{N}\to\mathbb{N}$ with the following
 properties:\/ $h_0(n)=o(n)$, and if\/ $h\ge h_0$ is any integer function,
 then every\/ $\Gamma$ of order\/ $n$ admits a zero-sum set\/ $A_0\subset \Gamma$
 such that\/ $|A_0|=h(n)$ and\/ $\Gamma\setminus A_0$ is partitionable
 into pairwise disjoint zero-sum subsets\/ $A_1, A_2,\ldots , A_t$ with\/
 $|A_i| = r_i$ whenever\/ $r_1 + r_2 + \ldots + r_t = n - |A_0|$ and\/
  {$r_i \geq 3$} for all\/ $1 \leq i \leq t$.
\end{cor}

Due to the possible strengthening indicated after Lemma \ref{l:FR},
 the above proof shows that only an overwhelming presence of values $r_i=3$
  can be responsible for the error term $\vep n$.
For this reason, on slightly restricted sequences of the $r_i$ we can obtain
 an almost optimal result.

\begin{cor}\label{cor:2}
 If the number of\/ $r_i=3$ is at most\/ $(1/3-c)\cdot n$ for a fixed\/ $c>0$,
  and the number of\/ $r_i=2$ does not exceed\/ $n/4$,
  then for sufficiently large\/ $n>n_c$ every sequence\/ $r_1,\dots,r_t$
  admits disjoint zero-sum subsets\/ $A_1, \ldots , A_t$ with\/
 $|A_i| = r_i$ in every\/ $\Gamma$ of order\/ $n\ge r_1 + \ldots + r_t +5$.
\end{cor}

 As in the preceding proof, $n/4$ can be replaced with $|R|/2$ also here.
On the one hand this condition depends on the actual $\Gamma$, while
 the restricted version given in the corollary is universally valid.
On the other hand the modified condition $|R|/2$ is best possible
 for every $\Gamma$.

\begin{rem}
 The bound on the number of pairs\/ $r_i=2$ in Theorem \ref{t:pairs2}
  is tight, because\/ $|R|=n/2$
  may occur, and then only\/ $n/4$ zero-sum pairs exist in\/ $\Gamma$.
\end{rem}

The following conjecture was raised recently.

\begin{conj}[\cite{ref_CicZ}]   \label{c:part}
 Let\/ $\Gamma$ of order\/ $n$ have  more than one involution.
For every partition\/ $n-1 = r_1 + r_2 + \ldots + r_t$ of\/ $n-1$ with\/
 {$r_i \geq 3$} for\/ $1 \leq i \leq t$ and for any possible positive integer\/ $t$,
  there is a partition of\/ $\Gamma-\{0\}$ into pairwise disjoint subsets\/
 $A_1, A_2,\ldots , A_t$ such that\/ $|A_i| = r_i$ and\/ $\sum_{a\in A_i}a = 0$
  for\/ $1 \leq i \leq t$.
\end{conj}

Note that the conjecture is true for $\Gamma\cong (\mathbb{Z}_2)^m$ as Egawa proved in \cite{Egawa}.
Moreover since for every $\Gamma$ having more than one involution we have
 $\sum_{g\in \Gamma}g=0$,
% (see e.g.\ \cite{CN}),
% {\komm (~\footnote{ZS -- this equality looks very basic; is \cite{CN}
% really the proper reference, as late as in 2004?} )}
%{\komm (~\footnote{SC -- yes it is a basic property and it's know for algebraists.
% I didn't find an earlier reference, but if you wish we can delete the reference
% and just leave it as known fact } )}
  the following two observations are valid by Theorem~\ref{glowne} and Corollary~\ref{cor:2}, respectively.

\begin{obs}\label{duze}Let\/ $\Gamma$ of order\/ $n$ have more than one involution. For
every partition\/ $n-1 = r_1 + r_2 + \ldots + r_t$ of\/ $n-1$, with\/ {$r_i \geq 3$} for\/ $1 \leq i \leq t$ and\/ $r_t\geq 3n/4$, for any possible positive integer\/ $t$, there is a partition of\/ $\Gamma-\{0\}$ into pairwise disjoint subsets\/ $A_1, A_2,\ldots , A_t$, such that\/ $|A_i| = r_i$ and\/ $\sum_{a\in A_i}a = 0$ for $1 \leq i\leq t$.
\end{obs}

%\begin{thm}\label{duze}
%Let $\Gamma$ be of order $n$, and consider any integers $r_1,r_2,\dots,r_t$ with
 %$n > 4 (r_1 + r_2 + \ldots + r_t)$ and {$r_i \geq 3$} for all $1 \leq i \leq t$.
%Then there exist pairwise disjoint subsets $A_1, A_2,\ldots , A_t$ in $\Gamma-\{0\}$
% such that $|A_i| = r_i$ and $\sum_{a\in A_i}a = 0$ for all $1 \leq i \leq t$.
%\end{thm}

%\textit{Proof.}   Therefore by Theorems~\ref{zerosum} and \ref{glowne} we are done.~\qed\bsk

\begin{obs}\label{duzen}Let\/ $\Gamma$ of large enough order\/ $n$ have more than one involution. For
every partition\/ $n-1 = r_1 + r_2 + \ldots + r_t$ of\/ $n-1$, with\/ {$r_i \geq 4$} for $1 \leq i \leq t$, for any possible positive integer\/ $t$, there is a partition of\/ $\Gamma-\{0\}$ into pairwise disjoint subsets\/ $A_1, A_2,\ldots , A_t$, such that\/ $|A_i| = r_i$ and\/ $\sum_{a\in A_i}a = 0$ for\/ $1 \leq i\leq t$.
\end{obs}

\textit{Proof.} If $n_t\geq 5$, we apply Corollary~\ref{cor:2}
 for $n_1,\dots,n_{t-1}$ to create the first $t-1$ sets.
The remaining $n_t$ elements of $\Gamma$ automatically sum up to zero,
 serving for the largest set.
Otherwise, if all $n_1=\ldots=n_t=4$, using the notation in the proof of
 Theorem \ref{t:pairs2} we partition $I\cup\{0\}\cong (\zet_2)^m$
 into zero-sum quadruples by Egawa's theorem, and partition $R$ into
 quadruples of the type $(a,b,-a,-b)$.~\qed

\section{Some applications}

Consider a simple graph $G=(V,E)$ whose order we denote by $n=|V|$.
% We denote by $V (G)$  the vertex set and $E(G)$  the edge set of $G$.
The \emph{open neighborhood} $N(x)$ of a vertex $x$ is the set of vertices adjacent to $x$, and the degree $d(x)$ of $x$ is $|N(x)|$, the order of the neighborhood of $x$. In this paper we  also investigate group  distance magic labelings, which belong to a large family of magic-type labelings. Generally speaking, a magic-type labeling of $G=(V,E)$ is a mapping from $V,$ $E,$ or $V\cup E$ to a set of labels which most often is a set of integers or group elements. The magic labeling (in the classical point of view)  with labels being the elements of an Abelian group has been studied for a long time (see papers by  Stanley~\cite{ref_Sta,ref_Sta2}).  Froncek in \cite{Fro1} defined the notion of group distance magic graphs, which are the graphs allowing a bijective labeling of vertices
with elements of an Abelian group resulting in a constant sum of neighbor labels.

A $\Gamma$\emph{-distance magic labeling} of a graph $G = (V, E)$ with $|V| = n$ is a bijection $\ell$ from $V$ to an Abelian group $\Gamma$ of order $n$
such that the weight $w(x) =\sum_{y\in N(x)}\ell(y)$ of every vertex $x \in V$ is equal to the same element $\mu\in \Gamma$, called the \emph{magic
constant}.

Notice that the constant sum partitions of a group $\Gamma$ lead to complete multipartite $\Gamma$-distance magic labeled graphs. For
instance, the partition $\{0\}$, $\{1, 2, 4\}$, $\{3, 5, 6\}$ of the group $\zet_7$ with constant sum $0$ leads to a $\zet_7$-distance magic labeling
of the complete tripartite graph $K_{1,3,3}$ (see \cite{ref_CicZ}). Using Theorem~\ref{duze} we are able to prove the following.

\begin{obs}\label{even} Let\/ $G = K_{n_1,n_2,\ldots,n_t}$ be a complete\/
 $t$-partite graph such that\/ $3\leq n_1\leq n_2 \leq\ldots\leq n_t$ and\/
  $n = n_1 +n_2 +\ldots+n_{t}$.
 Let\/ $\Gamma$ be an Abelian group of order\/ $n$ having more than three involutions.
The graph\/ $G$ is\/ $\Gamma$-distance magic whenever\/ $n_t\geq 3n/4-1$.
\end{obs}
\textit{Proof.}  There exists a zero-sum partition $A_1',A_2',\ldots,A_t'$  of the set
 $\Gamma- \{0\}$ such that $|A_t'|=n_t-1$ and $|A_i'| =n_i$ for every
  $1 \leq i \leq t-1$ by Theorem~\ref{duze}.
 Set $A_t=A_t'\cup \{0\}$  and $A_i=A_i'$ for every $1 \leq i \leq t-1$.
Label now the vertices from $V_i$, where $V_i$ is the vertex class of
 cardinality $n_i$, using the elements from the set $A_i$
 for $i\in\{1,2,\ldots,t\}$.~\qed\bsk

Analogously, for $n$ large enough, by Observation~\ref{duzen} we  can obtain:

\begin{obs}\label{even-large} Let\/ $\Gamma$ be an Abelian group of large enough
 order\/ $n$ having more than one involution.
If\/ $G = K_{n_1,n_2,\ldots,n_t}$ is a complete\/ $t$-partite graph such that\/
 $4\leq n_1\leq n_2 \leq\ldots\leq n_t$ and\/ $n = n_1 +n_2 +\ldots+n_{t}$,
    then\/ $G$ is a\/ $\Gamma$-distance magic graph.
\end{obs}
%\textit{Proof.}  If $n_t\geq 5$, we apply Corollary~\ref{cor:2}
% for $n_1,\dots,n_{t-1}$ to label the first $t-1$ vertex classes.
%The remaining $n_t$ elements of $\Gamma$ automatically sum up to zero,
% serving for the largest vertex class.
%Otherwise, if all $n_1=\ldots=n_t=4$, using the notation in the proof of
% Theorem \ref{t:pairs2} we partition $I\cup\{0\}\cong (\zet_2)^m$
 %into zero-sum quadruples by Egawa's theorem, and partition $R$ into
 %quadruples of the type $(a,b,-a,-b)$.\qed

\bibliographystyle{plain}

\end{document}